\documentclass{amsart} %

\usepackage{amsfonts,amssymb,amsmath,graphicx}
\usepackage[mathscr]{eucal}

\newtheorem*{thm*}{Theorem}
\newtheorem*{prop*}{Proposition}
\newtheorem*{lem*}{Lemma}

\def\overstrike#1#2{{\setbox0\hbox{$#2$}\hbox to \wd0{\hss
    $#1$\hss}\kern-\wd0\box0}}
\def\widedoublestrike#1{{\setbox0\hbox{$#1$}\hbox to \wd0{\hss
    $#1$\hss}\kern-0.67\wd0\box0}}
     \def\connsum{\raise.25ex\hbox{\overstrike\parallel=}}
     \font\bxtwelve=cmbx12 
     \def\plumb#1{\thinspace{{\lower.75ex%
          \hbox{\text{\bxtwelve*}}}}_{#1}\thinspace} 
     \def\Bd{\partial} 
     
     \def\sub{\subset} 
     \def\sgn{\operatorname{sgn}}
     \def\Z{\mathbb Z}

\begin{document}

\title
{A non-ribbon plumbing of fibered ribbon knots}
\author{Lee Rudolph}
\address{Department of Mathematics, Clark University, 
Worcester MA 01610 USA}  
\email{lrudolph@black.clarku.edu}


\begin{abstract} 
A closer look at 
an example introduced by Livingston \& Melvin 
and later studied by Miyazaki 
shows that a plumbing of two fibered ribbon knots 
(along their fiber 
surfaces) may be algebraically slice
yet not ribbon.
\end{abstract}

\maketitle

Trivially, the connected sum (i.e., $2$-gonal Murasugi sum) 
of ribbon knots is ribbon.  
Non-trivially \cite{Stallings,Gabai:Murasugi1},
any Murasugi sum of fibered knots (along their fiber surfaces) 
is fibered.  In light of these facts, perhaps the following 
is a bit surprising.

\begin{thm*} There exist fibered ribbon knots $K_1, K_2$ 
and an algebraically slice plumbing 
\emph{(\emph{$4$-gonal Murasugi sum})} $K_1\plumb{} K_2$, 
along fiber surfaces, which is not ribbon.
\end{thm*}

The proof uses a slight embellishment of a result from
\cite{Neumann-Rudolph}.
Following \cite{Neumann-Rudolph},
for any knot $K$, and relatively prime integers $m, n$ 
with $m\ge 1$, let $K\{m,n\}$ denote any simple closed 
curve on the boundary $\Bd N(K)$ of a tubular neighborhood $N(K)$ 
of $K$ in $S^3$ such that $K\{m,n\}$ represents $m$ times the
class of $K$ in $H_1(N(K);\Z)$ and has linking number $n$ with
$K$.  For instance, if $O$ denotes an unknot, then $O\{m,n\}$
is the (fibered) torus knot of type $(m,n)$.  
Abbreviate $K\{m,n\}\{p,q\}$ to $K\{m,n;p,q\}$.
(N.B.: there is no universally accepted standard 
notation for such iterated cable knots.  In particular,
instead of $O\{m,n;p,q\}$, Livingston \& Melvin 
\cite{Livingston-Melvin} write $(q,p;n,m)$,
and Miyazaki \cite{Miyazaki} writes $(p,q;m,n)$.) 

\begin{figure}
\begin{center}
\includegraphics[width=5in]{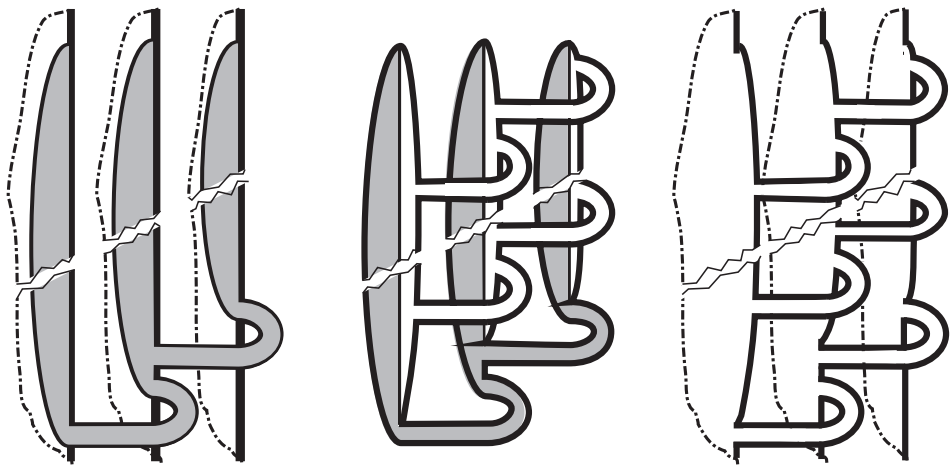}
\end{center}
\caption{A $3$-patch $N_1$ on $F_{(3,1)}$;
a $3$-patch $N_2$ on the (braided) fiber 
surface $D_{(3,n)}$; $F_{(3,n)}$ as a 
$6$-gonal Murasugi sum $F_{(3,1)}\plumb{h}D_{(3,n)}$
along the obvious diffeomorphism $h:N_1\to N_2$.}
\label{deplumbing a cable (m=3)}
\end{figure}

\begin{prop*} For any $K$, 
$K\{m,n\}=K\{m,\sgn n\}\plumb{} O\{m,n\}$ 
is a $2m$-gonal Murasugi sum, along a suitable
Seifert surface $F_{(m,\sgn n)}$ for $K\{m,\sgn n\}$ 
and a fiber surface $D_{(m,n)}$ for $O\{m,n\}$;
for fibered $K$, $F_{(m,\sgn n)}$ is a fiber surface.  
\end{prop*}
\begin{proof}
The case $m=1$ is trivial.  Take $m>1$, $\sgn n={\pm}1$.
Let $F\sub S^3$ be a Seifert surface for $K$.  Construct 
a Seifert surface $F_{(m,{\pm}1)}$ from $m$ parallel 
copies of $F$, each adjacent pair of copies joined,
in order, by a $1$-handle with a half-turn of sign $\pm$.
If $F$ is a fiber surface, then so is $F_{(m,{\pm}1)}$
(see \cite{Stallings}); in any case 
$\Bd F_{(m,\pm1)}=K\{m,{\pm}1\}$.  
The proof is finished by contemplating an appropriate figure 
(see Figure~\ref{deplumbing a cable (m=3)}).
\end{proof}

\begin{proof}[Proof of theorem.]
Livingston \& Melvin \cite{Livingston-Melvin} drew attention
to the connected sum of iterated torus knots
$K:=O\{2,3;2,13\}\connsum O\{2,15\}\connsum%
O\{2,-3;2,-15\}\connsum O\{2,-13\}$.  (Their
motivation was a vague question \cite{Rudolph:query}
about relations among the concordance classes of the 
knots associated to complex plane curve singularities,
and they gave an answer to one form of the question
by observing that $K$ is algebraically slice.)
Miyazaki \cite{Miyazaki} showed that $K$ is not ribbon.  
By the proposition, and the indifference of connected sums (of knots)
to the location of the summation,
\begin{equation*}\begin{split}
K & = ( O\{2,3;2,1\} \plumb{} O\{2,13\} ) \connsum O\{2,15\} 
      \phantom{O\{2,-3;2,-1\}00000001}\\
  & \qquad\qquad 
        \connsum ( O\{2,-3;2,-1\} \plumb{} O\{2,-15\} ) 
                                            \connsum O\{2,-13\} \\
\end{split}\end{equation*}

\begin{equation*}\begin{split}
  & = \bigl(( O\{2,3;2,1\} \connsum O\{2,-3;2,-1\})
      \plumb{} (O\{2,13\}\connsum O\{2,-13\})\bigr) \\
  & \qquad\qquad \plumb{} (O\{2,15\}\connsum O\{2,-15\})
\end{split}\end{equation*}
(where all the Murasugi sums are $4$-gonal and along fiber 
surfaces).  The connected sum of a knot and its mirror image 
is ribbon, so $R_1:=O\{2,3;2,1\} \connsum O\{2,-3;2,-1\}$,
$R_2:=O\{2,13\}\connsum O\{2,-13\}$, and 
$R_3$ $:=O\{2,15\}\connsum O\{2,-15\}$ are ribbon.
If $R_1\plumb{} R_2$ is ribbon, let $K_1:=R_1\plumb{} R_2$,
$K_2:=R_3$; if not, let $K_1:=R_1$, $K_2:=R_2$.
In either case, $K_1$ and $K_2$ are ribbon but $K_1\plumb{} K_2$
is not ribbon.  On the other hand, since (as mentioned 
in \cite{Neumann-Rudolph}) the Murasugi sum in the proposition 
is obviously ``direct'' (i.e., the Seifert linking between any cycle 
on $F_{(m,\sgn n)}$ and any cycle on $D_{(m,n)}$, taken in either 
order, is $0$), $K_1\plumb{} K_2$ is algebraically slice, for 
clearly a Murasugi direct sum is algebraically concordant to 
the connected sum of the same summands.
\end{proof}

\providecommand{\bysame}{\leavevmode\hbox to3em{\hrulefill}\thinspace}

\end{document}